\documentclass[reqno]{amsart}
\usepackage{amssymb,latexsym}

\newtheorem{thm}{Theorem}[section]
\newtheorem{lemma}{Lemma}[section]

\newtheorem{prop}{Proposition}[section]

\theoremstyle{definition}
\newtheorem{defn}{Definition}[section]
\theoremstyle{remark}
\newtheorem{rem}{Remark}[section]

\makeatletter
\@addtoreset{figure}{section}
\def\thefigure{\thesection.\@arabic\c@figure}
\def\fps@figure{h, t}
\@addtoreset{table}{bsection}
\def\thetable{\thesection.\@arabic\c@table}
\def\fps@table{h, t}
\@addtoreset{equation}{section}

\makeatother

\begin{document}

\title[Reduction and Covariant Euler-Poincar\'{e}]
{Reduction in Principal Fiber Bundles: Covariant Euler-Poincar\'{e} Equations}

\author[M. Castrill\'{o}n L\'{o}pez]{Marco Castrill\'{o}n L\'{o}pez}
\address{ Departmento de Geometr\'{i}a y Topolog\'{i}a, Universidad Complutense
 de Madrid, 28040 Madrid, Spain}
\email{mcastri@mat.ucm.es} 

\author[T.S. Ratiu]{Tudor S. Ratiu}
\address{Departement de Mathematiques\\ Ecole Polytechnique federale Lausanne\\
CH - 1015 Lausanne, Switzerland}
\email{ratiu@ratiu@masg1.epfl.ch}

\author[S. Shkoller]{Steve Shkoller}
\address{ {CNLS,MS-B258\\Los Alamos, NM 87545}}
\address{ {CDS\\California Institute of Technology, 107-81\\Pasadena,CA 91125}}
\email{shkoller@cds.caltech.edu}

\date{January 1998; this version, October 9, 1998. To appear in 
Proc. Amer. Math. Soc.}

\begin{abstract}
Let $\pi:P\rightarrow M^n$ be a principal $G$-bundle, and let ${\mathcal{L}}:J^1P \rightarrow
\Lambda^n(M)$ be a $G$-invariant Lagrangian density.   We obtain the Euler-Poincar\'{e}
equations for the reduced Lagrangian $l$ defined on ${\mathcal C}(P)$, the bundle of connections
on $P$.
\end{abstract}

\maketitle

%\tableofcontents

\section{Introduction}

Classical Euler-Poincar\'{e} equations arise through a reduction of the
variational principal
$\int_a^b L(\dot g(t))dt$ where $L:TG \rightarrow {\mathbb R}$ is a $G$-invariant
Lagrangian defined on the tangent bundle of a Lie group $G$.  In this setting,
one defines the reduced Lagrangian $l: TG/G \cong {\mathfrak g} \rightarrow 
{\mathbb R}$ by $l(\xi) = L(R_{g^{-1}} \dot g)$  (or by left-translation depending
on the Lagrangian), and proves that with a restricted class of variations, the
extremal $\xi$ of $\int_a^b l(\xi(t))dt$ is equivalent to the extremal of the
original variational problem for $L$.

The purpose of this note is to extend the variational reduction program to the
setting of a principle fiber bundle $\pi: P \rightarrow M$, using the fact 
that $J^1P/G \cong {\mathcal C}(P)$, where $J^1P$ is the first jet bundle of $P$, and
${\mathcal C}(P)$ denotes the bundle of connections on $P$. The reduced equations 
obtained can be seen as generalized Euler-Poincar\'{e} equations for field theory.
A remarkable fact is that these reduced equations on ${\mathcal C}(P)$ are not enough
for the reconstruction of the original problem for dim$M>1$. In classical mechanics, 
direct integration of the Euler-Poincar\'{e} equations gives solutions of the 
variational problem, but for field theory a set of compatibility equations are 
needed and they arise as the vanishing of the curvature of the reduced solution.
 This paper is the first in a series.  Herein, we establish the covariant reduction 
process in the case that $G$ is a matrix group.  In following notes, we shall make
the extension to more general Lie groups, as well as to the very interesting setting 
of homogeneous spaces.

\section{Preliminaries and notations}
Throughout this paper, differentiable will mean $C^{\infty}$ and if $E \to M$ is a fiber bundle,
$C^{\infty}(E)$ will denote the space of differentiable sections of $E$ over $M$.

\subsection{The bundle of connections}

Let $\pi \colon P \rightarrow M$ be a principal $G$-bundle. The right group action of $G$
on $TP$ is given by the lifted action
$$ X\cdot g=(R_{g})_{*}(X), \ \  \forall X\in TP, \ g \in G.$$ 
The quotient $TP/G$ is a differentiable manifold
and is endowed with a vector bundle structure over $M$. Let $\text{ad}P:=
(P\times {\mathfrak g})/G$, the bundle associated to $P$ by the adjoint
representation of $G$ on ${\mathfrak g}$.  With $VP$ the vertical
subbundle of $TP$, 
the map $h\colon \text{ad}P \rightarrow VP/G$  given by
$$ h((p,\xi )_{G})=(\hat{\xi}_{p})_{G} $$
is a vector bundle diffeomorphism, where $\hat{\xi}_p=
(d/dt)|_0 p\cdot \exp (t\xi )$. 
Let 
\begin{equation}
\ddagger :VP\to \mathrm{ad}P\stackrel{h}{\cong }VP/G  \label{1}
\end{equation}
be the projection induced by the diffeomorphism $h$. The fibers $(\mathrm{ad}%
P)_{x}$ of the adjoint bundle are endowed with a Lie algebra structure
determined by the following condition 
\begin{equation}
\lbrack (p,\xi )_{G},(p,\eta )_{G}]=(p,[\xi ,\eta ])_{G},\hspace{5mm}\forall
p\in \pi ^{-1}(x),\forall \xi ,\eta \in \mathfrak{g},  \label{2}
\end{equation}
where $[\cdot,\cdot]$ denotes the bracket on ${\mathfrak g}$.

The quotient modulo $G$ of the following exact sequence of vector
bundles over $P$, 
$$ 0\to VP\to TP\stackrel{\pi _{*}}{\longrightarrow }\pi ^{*}TM\to 0,$$
becomes the exact sequence of vector bundles over M 
$$ 0\to \text{ad}P\to TP/G\stackrel{\pi _{*}}{\longrightarrow }TM\to 0,$$
which is called the \emph{Atiyah sequence\/} (see, for example \cite{Atiyah}).

\begin{defn}
A connection on $P$ is a distribution ${\mathcal H}$ complementary to $VP$,
such that ${\pi_*}_p: {\mathcal H}_p \rightarrow T_{\pi(p)}M$ is an
isomorphism for all $p \in P$.  The horizontal lift of a vector field $X$ on
$M$ is the vector field $\tilde{X}$ on $P$ defined by
$\tilde{X}(p) := (\pi_*|_{{\mathcal H}_p})^{-1} X(\pi(p))$.
\end{defn}

Let ${\mathcal H}$ be a connection on $P$, and let $\tilde{X}\in {\mathfrak X}(P)$ be
the horizontal lift with respect to ${\mathcal H}$ of a vector field $X\in \mathfrak{X}(M)$.
The horizontal lift is a $G$-invariant vector field on $P$ projecting onto 
$X$. Namely, 
\begin{equation}\label{s1}
{R_g}_* {\mathcal H}_p = {\mathcal H}_{pg} \ \  \forall g \in G \text{ and }
p \in P.
\end{equation}
Hence, there exists a splitting of the Atiyah sequence 
$$ \sigma \colon TM \rightarrow TP/G,\hspace{5mm}\sigma (X)=\tilde{X}.$$
Conversely, any splitting $\sigma \colon TM \rightarrow TP/G$ induces a unique
connection:  let $\xi \in {\mathfrak g}$ and $\psi \in {\mathcal H}$, and
define the ${\mathfrak g}$-valued $1$-form ${\mathcal A}$ on $P$ by
$$ {\mathcal A} \langle \hat{\xi} + \psi \rangle = \xi.$$
It follows that
$$ {\mathcal H}_p = \text{Ker}{\mathcal A}_p,$$
so there is a natural bijective correspondence
between connections on $P$ and splittings of the Atiyah sequence.

In the case that $P=M \times G$ is trivial, condition (\ref{s1}) implies
that the horizontal lift of a vector field $X$ on $M$ is given by
\begin{equation}\label{s2}
\tilde{X}=(X_x, - ({\mathcal A}_x\langle X \rangle g)_g)=
(X_x, -(R_g)_* {\mathcal A}_x \langle X \rangle), \ \ \ p=(x,g).
\end{equation}

\begin{defn}
We denote by $p:{\mathcal C}(P) \rightarrow M$
the subbundle of $\text{Hom}(TM,TP/G)$
determined by all linear mappings  (see \cite{Eck}, \cite{Garcia})
$$\sigma_x: T_xM \rightarrow (TP/G)_x \text{ such that }
\pi_*\circ \sigma_x=\text{Id}_{T_xM}.$$
\end{defn}

An element $\sigma_x\in {\mathcal C}(P)_x$
is  a distribution at $x$; that is, $\sigma _{x}$
induces a complementary subspace ${\mathcal H}_p$ of the vertical subspace
$V_pP$ for any $p\in \pi^{-1}(x)$. 
Addition of a linear mapping $l_{x}\colon T_{x}M\to (\mathrm{ad}P)_{x} \in 
\text{Ker}\pi_*$ to $\sigma _{x}$, produces
another element $\sigma _{x}^{\prime }=l_{x}+\sigma _{x}\in (\mathcal{%
C}(P))_{x}$, so that ${\mathcal C}(P)$ is an affine bundle modeled over 
the vector bundle 
$\mathrm{Hom}(TM,\text{ad}P)\simeq T^{*}M\otimes \mathrm{ad}P$.

Accordingly, any global section $\sigma \in C^\infty(\mathcal{C}(P))$ can be
identified with a global connection on $P$.  Similarly, the difference of the two 
global sections $\sigma$ and ${\mathcal H}$ may be identified with a section of 
the bundle $T^*M\otimes \text{ad}P$. If we fix a connection ${\mathcal H}$, the map 
\begin{equation}
\Phi _{{\mathcal H}}: {\mathcal C}(P) \rightarrow T^*M\otimes \text{ad}P
\text{ given by }
\Phi _{{\mathcal H}}(\sigma)=\sigma -{\mathcal H},  \label{3}
\end{equation}
is a fibered diffeomorphism. Note, however, that although $\mathcal{C}%
(P)\simeq T^{*}M\otimes \mathrm{ad}P$, the diffeomorphism is not canonical;
it depends on the choice of the connection ${\mathcal H}$. We will denote by 
$\sigma^{\mathcal H}$ the image of $\sigma$ under $\Phi_{\mathcal H}$.

\subsection{The identification $J^1P/G \simeq \mathcal{C}(P)$}

\begin{defn}
Let $\pi: P\rightarrow M$ be a principal $G$-bundle and denote the
$1$-jet bundle of local section of $\pi$  by $\pi _{1}\colon
J^{1}P\to M$. This is the affine bundle of all linear mappings
$\lambda_x:T_xM \rightarrow T_pP$ such that ${\pi _*}_p\circ \lambda_x=
\text{Id}_{T_xM}$ for any $p \in \pi^{-1}(x)$.  If $s$ is a local section
of $P$, its first jet extension $j^1s$ is identified with the tangent
map of $s$, i.e. $j^1_x s = T_x s$, $x \in M$.
\end{defn}
The group $G$ acts on $J^{1}P$ in a natural way by 
\begin{equation}
j_{x}^{1}s\cdot g=j_{x}^{1}(R_{g}\circ s),  \label{4}
\end{equation}
where $R_{g}$ is the right action of $G$ on $P$. 
The quotient $J^1P/G$
exists as a differentiable manifold and can be identified with the bundle of
connections in the following way. We have $j_{x}^{1}s \cdot g=j_{x}^1(R_g \circ s)
=T_x(R_g \circ s)=(R_g)_*T_xs$; then a coset $(j_{x}^{1}s)_G \in J^1P/G$ can be seen
as a $G$-invariant horizontal distribution over $M$, that is, an element in  
$\mathcal{C}(P)_{x}$.
Let 
\begin{equation}
q:J^{1}P\to \mathcal{C}(P)\simeq J^{1}(P)/G  \label{5}
\end{equation}
be the projection.
Let $U\subset M$ be a local neighborhood of $x \in M$.
If $s \in C^\infty(P|_U)$,
we obtain a local section $\sigma: U \rightarrow {\mathcal C}(P)$ as
$\sigma (x)=q(j_{x}^{1}s)$.

\section{Euler-Poincar\'{e} reduction}

Let $\pi :P\to M$ be a fiber bundle. A \emph{Lagrangian density} is a bundle
map $\mathcal{L}\colon J^{1}P\to \Lambda ^{n}M$, 
where $n=\dim M$. 
\begin{defn} A variation of $s\in C^\infty(P)$ is a curve
$s_\epsilon  = \phi_\epsilon \circ s$,  where $\phi_\epsilon$ is the flow
of a  vertical vector field $V$ on $P$ which is compactly supported in $M$.
One says that $s$ is a fixed point of the variational problem associated with
${\mathcal L}$ if
\begin{equation}
\delta \int_M{\mathcal L}(j^1s) :=
\left.\frac{d}{d\epsilon} \left[ \int_M {\mathcal L}(j^1s_\epsilon) \right]
\right|_{\epsilon=0} = 0  \label{6}
\end{equation}
for all variations $s_\epsilon$ of $s$.
\end{defn}

For a fixed volume form $dx$ on $M$, we define the \emph{Lagrangian}
associated to $\mathcal{L}$ as the mapping $L:J^{1}P\to {\mathbb R}$ which
verifies $\mathcal{L}(j_{x}^{1}s)=L(j_{x}^{1}s)dx$, $\forall j_{x}^{1}s\in
J^{1}P$. Then, formula (\ref{6}) becomes 
\[
\delta \int_{M}L(j_{x}^{1}s)dx=0. 
\]

Henceforth, we shall restrict attention to the principal $G$-bundle
$\pi :P\to M$ with $\dim M=n$ and with volume form $dx$.
\begin{defn}
A Lagrangian $L:J^{1}P\to {\mathbb R}$ is $G$-invariant if 
$$L(j_{x}^{1}s\cdot g)=L(j_{x}^{1}s),\hspace{5mm}\forall j_{x}^{1}s\in
J^{1}P,\forall g\in G,$$
where the action on $J^{1}P$ is defined in formula (\ref{4}). 
\end{defn}

If $L$ is a $G$-invariant Lagrangian, it defines a mapping 
$$ l:J^1P/G \simeq {\mathcal C}(P) \rightarrow {\mathbb R}$$
in a natural way.  With $\delta s := (d/d \epsilon)|_0  s_\epsilon \in
C^\infty(VP)$, define $\eta \in C^\infty(\text{ad}P)$ by
$$\eta(x) = \ddagger \delta s (x).$$

\begin{prop}\label{prop1}
Let $\pi :P \rightarrow M$ be a principal $G$-bundle, $G$ a matrix group,
with a fixed connection
${\mathcal H}$, and consider the curve $\epsilon \mapsto s_\epsilon = \phi_\epsilon \circ
s$, where $\phi_\epsilon$ is the flow of a $\pi$-vertical vector field $V$.
Define $\sigma_\epsilon =q(j^{1}s_\epsilon)$ and $\sigma^{\mathcal H}=\Phi
_{\mathcal H}(\sigma )$. Then
$$\delta \sigma :=(d/d \epsilon)|_0 \sigma_\epsilon
=\nabla^{\mathcal H}\eta -[\sigma^{\mathcal H}\langle \cdot \rangle,\eta ]$$
where $[\cdot, \cdot]$ is given by \textrm{(\ref{2})}, and
$\nabla^{\mathcal H}:C^\infty(\mathrm{ad}P) \rightarrow C^\infty(
T^*M \otimes \mathrm{ad}P)$
is the covariant derivative induced by ${\mathcal H}$ in the
associated bundle $\text{ad}P$ defined in a trivialization by
\begin{equation}\label{s5}
\nabla^{\mathcal H}\eta = T\xi + [{\mathcal A} \langle \cdot \rangle, \xi],
\end{equation}
where $\eta(x) = (x, \xi(x))$.
\end{prop}

%\begin{rem}
%Note that $\delta \sigma$ is a section of $T^*M \otimes \mathrm{ad}P$ as $%
%\sigma$ is a section of $\mathcal{C}(P)$ and the bundle of connections is an
%affine bundle modeled over the vector bundle $T^* M \otimes \mathrm{ad}P$.%
%\end{rem}
\begin{rem}
If ${\mathcal H}^{\prime}$ is another connection on $P$, then
$$ \nabla ^{{\mathcal H}^{\prime}} \eta -
[\sigma ^{{\mathcal H}^{\prime}}\langle \cdot \rangle,\eta ]
=\delta \sigma=\nabla ^{\mathcal H} \eta 
-[\sigma ^{\mathcal H}\langle \cdot \rangle ,\eta ]. $$
\end{rem}

\begin{rem}
If we consider a principal fiber bundle with a left group action instead of a
right action, then the expression for the infinitesimal variation is
$$ \delta \sigma =\nabla ^{{\mathcal H}}\eta +[\sigma ^{{\mathcal H}}\langle
\cdot \rangle , \eta]. $$
\end{rem}

\begin{proof}
Since this is a local statement, we may assume that
$P=U \times G$, where $U \subset M$
is open, with $\pi$ the projection onto the first factor, and
with right action $R_{g'}$ given by 
$$R_{g'} (x,g) = (x,g) \cdot g' = (x,g g').$$
Hence, $\text{ad}P \simeq M\times {\mathfrak g}$ via the map $((x,e),\xi
)_{G}\mapsto (x,\xi )$ and the projection $\ddagger :V_{(x,g)}P\to (\text{ad}P)_x
\simeq {\mathfrak g}$ is given explicitly by right translation
\begin{equation}\label{s3}
\ddagger (0_{x},v)=(R_{g^{-1}})_{*}v=vg^{-1},\hspace{5mm}\forall v\in
T_{g}G. 
\end{equation}

We identify the map $g \in C^\infty(U,G)$ with $s \in C^\infty(U \times M)$
by $s(x) = (x,g(x))$ and the map $\xi \in C^\infty(U,{\mathfrak g})$ with 
$\eta \in C^\infty(\text{ad}P)$ by $\eta(x) = (x,\xi(x))$. We have the
following identifications:
$$ (TP/G)_x\simeq T_{(x,e)}P\simeq T_xM\times T_eG\simeq T_xM\times 
{\mathfrak g},$$
so that
$$
\sigma_x=q(T_xs)=(\text{Id}_{T_{x}M},(R_{g^{-1}})_{*} T_xg)
=(\text{Id}_{T_xM},T_xg\cdot g^{-1}). $$
Then, 
\begin{eqnarray*}
\delta \sigma(x)&=&(d/d \epsilon)|_0 \sigma_\epsilon(x)=
(d/d \epsilon)|_0 (\text{Id}_{T_{x}M},T_{x}g_\epsilon \cdot g_{\epsilon }^{-1}) \\
&=&(0_{x},\{ (d/d \epsilon)|_0 T_xg_\epsilon \} \cdot
g^{-1}-T_{x}g\cdot g^{-1}\cdot \delta g\cdot g^{-1})\\
&=&(0_{x},T_x \delta g\cdot g^{-1}-T_xg\cdot g^{-1}\cdot \delta g\cdot
g^{-1})\\
&=&\left(0_{x},[\delta g\cdot g^{-1},T_{x}g\cdot g^{-1}\langle \cdot \rangle]
-\delta g\cdot g^{-1}\cdot T_{x}g\cdot g^{-1}+T_{x}\delta g\cdot g^{-1}\right)\\
&=&\left(0_{x},[\delta g\cdot g^{-1},T_{x}g\cdot g^{-1}\langle \cdot \rangle]
+T_{x}(\delta g\cdot g^{-1})\right),
\end{eqnarray*}
where the bracket is the commutator of matrices as $G$ is a matrix group.
Hence, $\delta \sigma (x)\in T_{x}^{*}M\otimes 
(\text{ad}P)_x\simeq T_{x}^{*}M \otimes {\mathfrak g}$, a 
${\mathfrak g}$-valued (vertical-valued) $1$-form.  Now, $\eta = 
\ddagger \delta s$, so using (\ref{s3}), $\xi = \delta g g^{-1}$. (We
make the identification $T_{\xi(x)}{\mathfrak g}\simeq{\mathfrak g}$.)

So for any vector field $X$ on $M$,
\begin{eqnarray*}
\delta \sigma \langle X \rangle &=& [\xi, \sigma\langle X \rangle] + T\xi \\
&=& [\xi, \sigma\langle X \rangle - \tilde{X}] + T\xi + [\xi,\tilde{X}].
\end{eqnarray*}
Let ${\mathcal A}$ be the local connection $1$-form associated to 
${\mathcal H}$. Then using (\ref{s2}),
$$ \delta \sigma = [\xi , \sigma\langle \cdot \rangle + {\mathcal A}\langle
\cdot \rangle] + T\xi - [\xi, {\mathcal A} \langle \cdot \rangle]
= -[\sigma^{\mathcal H}\langle \cdot \rangle, \xi] + T\xi - [\xi,
{\mathcal A}\langle \cdot \rangle ] .$$

Now to obtain the formula for $\nabla^{\mathcal H}$, we use the injective
correspondence between $C^\infty(\text{ad}P)$ and
$\{ f_\eta \in C^\infty(P,{\mathfrak g})| f_\eta(pg)= \text{Ad}_{g^{-1}}f_\eta(p),
\ p \in P, g\in G\}$. Hence,
$$ f_\eta(x,g) = \mathrm{Ad}_{g^{-1}} \xi(x).$$
It is standard (see \cite{KN}) that $\nabla^{\mathcal H} _X\eta$ is given
by $(f_\eta)_* \langle \tilde{X} \rangle$, so we need only use (\ref{s2}) to
compute the horizontal lift $\tilde{X}_{(x,e)}$.  We have that
\begin{eqnarray*}
(f_\eta)_* \langle \tilde{X} \rangle &=& (f_\eta)_* \langle X \rangle
+ (d/dt)|_0 f_\eta(x, \mathrm{exp}(-t {\mathcal A} \langle X \rangle) \\
&=&T\xi \langle X \rangle + (d/dt)|_0 \mathrm{Ad}_{\mathrm{exp}(t{\mathcal A}\langle
X \rangle)} \xi(x) = T\xi\langle X \rangle + [ {\mathcal A} \langle X \rangle , \xi],
\end{eqnarray*}
so that $\nabla^{\mathcal H}  \eta = T\xi + \text{ad}_{{\mathcal A}\langle\cdot\rangle} \xi$,
and this completes the proof.
\end{proof}

\begin{rem}
The dual of the adjoint bundle $({\mathrm{ad}}P)^*$ can be seen also as the bundle 
associated to $P$ by the dual adjoint representation of $G$ on ${\mathfrak g}^{*}$;
\emph{i.e.}, $g\mapsto {\mathrm{Ad}}_{g^{-1}}^{*}$. 
There is a similar injective correspondence between $C^\infty((\text{ad}P)^*)$ and the set
$\{ f_\eta \in C^\infty(P,{\mathfrak g})| f_\eta(pg)= \text{Ad}_g^* f_\eta(p),
\ p \in P, g\in G\}$. Hence, if $\nu \in C^\infty((\text{ad}P)^*)$ and
$\nu(x) = (x, \psi(x))$, then
the covariant derivative induced by the connection ${\mathcal H}$ in $(\text{ad}P)^*$ is defined by
\begin{equation}\label{s6}
\tilde{\nabla}^{\mathcal H}  \nu 
= T\psi - \text{ad}^*_{{\mathcal A}\langle\cdot\rangle} \psi,
\end{equation}
or equivalently
$$(\tilde{\nabla}_X \nu)\eta=X(\langle \nu,\eta\rangle)-\langle \nu,\nabla _X \eta
\rangle, \hspace{4mm}\forall \eta \in C^{\infty}({\mathrm{ad}}P), X\in {\mathfrak X}(M).
$$
\end{rem}

Given a section $\sigma \in \mathcal{C}(P)$, the mapping 
$l:\mathcal{C}(P) \rightarrow {\mathbb R}$ defines a linear operator 
$$\frac{\delta l}{\delta \sigma }:
T^{*}M\otimes \mathrm{ad}P \rightarrow {\mathbb R}$$
by 
$$\frac{\delta l}{\delta \sigma }(\zeta _{x})=\lim_{\epsilon \rightarrow 0}
\frac{l(\sigma(x)+\epsilon \zeta _{x})-l(\sigma (x))}{\epsilon },\hspace{5mm}%
\forall \zeta _{x}\in (T^{*}M\otimes \mathrm{ad}P)_{x}. $$
The operator $\delta l/\delta \sigma $ can be seen as a section of the dual
bundle $(T^{*}M\otimes \mathrm{ad}P)^{*}\simeq TM\otimes (\mathrm{ad}P)^{*}$%
. 

\begin{lemma} \label{lemma1}
For a fixed connection ${\mathcal H}$ on 
$\pi:P \rightarrow M$, there exists an associated divergence operator 
$\mathrm{div}^{\mathcal H}:C^\infty(TM \otimes (\mathrm{ad}P)^{*})\rightarrow 
C^\infty((\mathrm{ad}P)^{*})$ which satisfies the following conditions. Let
$\mathcal{X}, \mathcal{X}^{\prime }\in C^\infty(TM\otimes (\mathrm{ad}P)^*)$,
$\eta \in C^\infty(\mathrm{ad}P)$, and $f\in C^\infty(M)$.  Then
\begin{enumerate}
\item[i)]  $\mathrm{div}^{\mathcal H}(\mathcal{X+X}^{\prime })=
\mathrm{div}^{\mathcal H}(\mathcal{X})+\mathrm{div}^{\mathcal H}
(\mathcal{X})$, 

\item[ii)]  $\mathrm{div}^{\mathcal H}(f\mathcal{X})=\mathcal{X\cdot }\mathrm{d}f+f%
\mathrm{div}^{\mathcal H}(\mathcal{X})$,

\item[iii)]  $\mathrm{div}(\mathcal{X}\cdot \eta )=
(\mathrm{div}^{\mathcal H}\mathcal{X})\cdot \eta +
\mathcal{X}\cdot \nabla^{\mathcal H}\eta.$
\end{enumerate}
Furthermore, if $\{E^{1},...,E^{m}\}$ is a basis of local sections of the
bundle $(\mathrm{ad}P)^{*}$ for which any element $\mathcal{X}\in C^\infty(TM\otimes (%
\mathrm{ad}P)^{*})$ may be expressed as $\mathcal{X}=\sum X_{i}\otimes
E^{i}$, $X_{i}\in \mathfrak{X}(M)$, then
\begin{equation}
\mathrm{div}^{\mathcal H}(\mathcal{X})=\sum_{i=1}^{m}\left( \mathrm{div}%
(X_{i})\otimes E^{i}+\tilde{\nabla}_{X_{i}}^{\mathcal H}E^{i}\right) .  \label{8}
\end{equation}
\end{lemma}

\begin{rem}
In the case $P=M\times G$ and ${\mathcal H}$ is the trivial connection, then
$\mathrm{div}^{\mathcal H}$ is the usual divergence operator.
\end{rem}

\begin{proof}
We use the same notation as in the proof of Proposition \ref{prop1}.
Let $\{E_{1},\ldots ,E_{m}\}$ be a basis of $\mathfrak{g}$ and $\{E^{1},\ldots
,E^{m}\}$ its dual basis. Let $\mathcal{X}=\sum X_{i}\otimes E^{i}$ be any
section of $TM\otimes (\mathrm{ad}P)^{*}\cong TM\otimes \mathfrak{g}^{*}$, and 
let $\xi =\sum f^{i}\otimes E_{i}$ be any section of $\mathrm{ad}P\cong M\times 
\mathfrak{g}$. Using (\ref{s5}) and (\ref{s6}), we have that
\begin{eqnarray*}
\mathcal{X}\cdot \nabla^{\mathcal H}\eta &=&
\sum_{i=1}^{m}\left(Tf^i \langle X_i \rangle +
f^i\sum_{j=1}^{m} \langle  [\mathcal{A}\langle X_j\rangle,E_i],E^j \rangle\right)\\
&=&\sum_{i=1}^{m}\left( X_i(f^i)+f^i\sum_{j=1}^{m} 
\langle\mathrm{ad}_{\mathcal{A}\langle X_j \rangle}^* E^j,E_i \rangle\right)\\
&=&\sum_{i=1}^{m}\left( \mathrm{div}(f^iX_i)-f^i\mathrm{div}X_{i}\right)
+\sum_{j=1}^{m}\langle \mathrm{ad}_{\mathcal{A}\langle X_j\rangle}^*E^j,\eta \rangle\\
&=&\mathrm{div}(\mathcal{X\cdot }\eta )-\left( \langle\sum_{j=1}^{m}
\mathrm{div}X_j \otimes E^j,\eta \rangle + 
\langle\sum_{j=1}^{m}\tilde{\nabla}_{X_j}^{\mathcal H}E^j,\eta \rangle\right) ,
\end{eqnarray*}
where $\langle \cdot,\cdot \rangle$ is the natural pairing between 
$\mathfrak{g}$ and $\mathfrak{g}^*$.  

Hence, the operator $\mathrm{div}^{\mathcal H}$ satisfying iii) is
$$ \mathrm{div}^{\mathcal H}\left(\sum_{j=1}^{m}
X_j \otimes E^{j}\right)= \sum_{j=1}^{m}
(\mathrm{div}X_j\otimes E^i+\tilde{\nabla}_{X_j}^{\mathcal H}E^j).$$
This expression can be defined globally and it is straightforward to
verify items i) and ii).
\end{proof}

\subsection{Reduction}

\begin{thm}\label{thm1}
Let $\pi :P\to M$ be a principal $G$-fiber bundle over a manifold $M$ with a
volume form $dx$ and let $L:J^{1}P\to {\mathbb R}$ be a $G$ invariant
Lagrangian. Let $l:\mathcal{C}(P)\to {\mathbb R}$ be the mapping defined by $L$ in
the quotient. For a section $s:U\to P$ of $\pi $ defined in a neighborhood $%
U\subset P$, let $\sigma :U\to \mathcal{C}(P)$ be defined by $\sigma
(x)=q(j_{x}^{1}s)$. Then, for every connection ${\mathcal H}$ of the bundle 
$\pi |_{U}$, the following are equivalent:

\begin{itemize}
\item[1)]  $s$ satisfies the Euler-Lagrange equations for $L$,

\item[2)]  the variational principle
$$ \delta \int_{M}L(j_{x}^{1}s)dx=0$$
holds, for variations with compact support,

\item[3)]  the Euler-Poincar\'{e} equations hold: 
$$\mathrm{div}^{\mathcal H}\frac{\delta l}{\delta \sigma }=-\mathrm{ad}_{\sigma
^{\mathcal H} \langle \cdot \rangle}^*\frac{\delta l}{\delta \sigma },$$

\item[4)]  the variational principle 
$$ \delta \int_{M}l(\sigma (x))dx=0$$
holds, using variations of the form 
$$ \delta \sigma =\nabla ^{\mathcal H}\eta -
[\sigma^{\mathcal H}\langle \cdot \rangle,\eta ] $$
where $\eta :U\to \mathrm{ad}P$ is a section with compact support.

\end{itemize}
\end{thm}

\begin{proof}
1) $\Leftrightarrow $ 2) is a standard argument in the calculus of
variations.  For
2) $\Leftrightarrow $ 4), we use that
$$ \delta \int_{M}L(j_{x}^{1}s)dx=\delta \int_{M}l(\sigma (x))dx$$
with Proposition \ref{prop1}.

For 3) $\Leftrightarrow $ 4), we have that
$$ 0=\delta \int_{M}l(\sigma (x))dx=\int_{M}\frac{\delta l}{\delta \sigma }%
\delta \sigma dx=\int_{M}\frac{\delta l}{\delta \sigma }(\nabla ^{\mathcal H}\eta
-[\sigma^{\mathcal H}\langle \cdot \rangle,\eta ])dx. $$

Item iii) of Lemma \ref{lemma1} gives that
$$\frac{\delta l}{\delta \sigma }\nabla^{\mathcal H}\eta =\mathrm{div}(\frac{\delta l}
{\delta \sigma }\eta )-\mathrm{div}^{\mathcal H}(\frac{\delta l}{\delta \sigma })\eta,$$
so that
$$ 0=\int_{M}(\mathrm{div}(\frac{\delta l}{\delta \sigma }\eta )-\mathrm{div}%
^{\mathcal H}(\frac{\delta l}{\delta \sigma })\eta -\mathrm{ad}_{\sigma ^{\mathcal H}
 \langle \cdot \rangle}^{*}%
\frac{\delta l}{\delta \sigma }\eta )dx$$
As $\eta $ has compact support, by Stokes theorem, $\int_{M}\mathrm{div}
(\frac{\delta l}{\delta \sigma }\eta )dx=0$, so we conclude  that
$$ 0=\int_{M}(\mathrm{ad}_{\sigma ^{\mathcal H}\langle \cdot \rangle }^{*}
\frac{\delta l}{\delta \sigma }+%
\mathrm{div}^{\mathcal H}(\frac{\delta l}{\delta \sigma }))\eta dx,$$
for all sections $\eta$ of $\mathrm{ad}P$ with compact support. Thus, we
obtain the Euler-Poincar\'{e} equations.
\end{proof}

\begin{rem}
If we consider a principal fiber bundle with a left action, instead of a
right action, and a left invariant Lagrangian $L$, the Euler-Poincar\'{e}
equations are 
$$ \mathrm{div}^{\mathcal H}\frac{\delta l}{\delta \sigma}=\mathrm{ad}^{*}_{\sigma
^{\mathcal H}\langle \cdot \rangle}\frac {\delta l}{\delta \sigma}.$$

\end{rem}

\subsection{Reconstruction}

Let $s:U \to P$ be a solution of the variational problem defined by a $G$
invariant Lagrangian $L$. Then, the section $\sigma =q(j^1 s)$ of the bundle
of connections is a solution of the Euler-Poincar\'{e} equations (Theorem
\ref{thm2}). This new section is a connection which verifies that $s(U)$ is an
integral manifold, that is, $\sigma $ is a flat connection. Conversely,
given a flat connection $\sigma $ which verifies the Euler-Poincar\'{e}
equations, the integral submanifolds in $P$ of $\sigma $ are the image of
the sections of the solution of the original variational problem. In other words

\begin{thm}\label{thm2}
The following systems of equations are equivalent:

\begin{enumerate}
\item[i)]  Euler-Lagrange equations of $L$, and

\item[ii)]  Euler-Poincar\'{e} equations of $l$ together with vanishing
curvature.
\end{enumerate}

The projection of a solution $s$ of $\mathrm{i)}$ gives a solution $\sigma
=q(j^{1}\sigma )$ of $\mathrm{ii)}$, and the integral manifolds of a solution $\sigma $
of $\mathrm{ii)}$ provides a solution of $\mathrm{i)}$.
\end{thm}

That is, the Euler-Poincar\'{e} equations are not sufficient for reconstructing
the solution of the original variational problem. One must impose an
additional compatibility condition given by the vanishing of the curvature.  See
\cite{MPS} for additional discussion.

\section{Examples of reduction in a principal fiber bundle.}

\subsection{Classical Euler-Poincar\'{e} equations.} 
For a Lie group $G$, we consider the principal fiber bundle $\pi:{\mathbb R}
\times G \rightarrow {\mathbb R}$, where $\pi$ is the projection onto the
first factor.
Let $L:J^{1}P\simeq {\mathbb R}\times TG\to {\mathbb R}$ be a right $G$%
-invariant Lagrangian. We fix the trivial connection and obtain the
following identifications: 
$$ \mathcal{C}(P)\simeq T^{*}{\mathbb R}\otimes \mathrm{ad}P\simeq 
({\mathbb R}\times ){\mathbb R}\otimes ({\mathbb R}\times \mathfrak{g})
\simeq {\mathbb R}\times \mathfrak{g}.$$

Again, we identify
$s\in C^\infty(P)$, $\eta \in C^\infty(\mathrm{ad}P)$ and 
$\sigma \in C^\infty (\mathcal{C}(P))$  with the maps
$g \in C^\infty({\mathbb R},G)$, 
$\eta \in C^\infty({\mathbb R},{\mathfrak g})$, and $\sigma
\in C^\infty({\mathbb R},{\mathfrak g})$, respectively.
Because of the trivial connection, $\mathrm{div}^{\mathcal H}$ is simply
the usual divergence operator satisfying
$\mathrm{div}(f\frac{\partial }{\partial t})=\frac{df}{dt}$.

We recover the classical the Euler-Poincar\'{e} equations
$$ \frac{d}{dt}\frac{\delta l}{\delta \sigma }=
-\mathrm{ad}_{\sigma }^{*} \frac{\delta l}{\delta \sigma }$$
for a right invariant Lagrangian (see \cite{MR}).

\subsection{Harmonic maps.}
Let $(M,g)$ be a compact oriented $C^\infty$
$n$ dimensional Riemannian manifold, and let $(G,h)$ 
be an $m$ dimensional Riemannian  matrix Lie group. With $P=M\times G$, we
denote the principal fiber bundle by $\pi:P \rightarrow M$, and by
triviality, identify $C^\infty(P)$ with $C^\infty(M,G)$.  For each
$\phi \in C^\infty(M,G)$, the Riemannian
metrics on $M$ and $G$ naturally induce a metric $\langle \cdot, \cdot \rangle$
on $C^\infty(T^*M \otimes \phi^*(TG))$, and so we may define the energy 
${\mathcal E}$ on $C^\infty(M,G)$ by
\begin{equation}\label{t1}
{\mathcal E}(\phi) = \int_M L(j^1\phi)dx, \text{  where  }
L(j^1\phi)={\frac{1}{2}} \langle T\phi, T\phi \rangle.
\end{equation}
The Euler-Lagrange equations for (\ref{t1}) are given by
\begin{equation}\label{t2}
{\mathrm {Tr}}(\nabla T\phi )=0, 
\end{equation}
where $\nabla$ is the induced Riemannian covariant derivative on
$C^\infty(T^*M \otimes \phi^*(TG))$ and Tr is the trace defined by $g$ (see, for example \cite{EL}). By definition, the set of harmonic 
maps from $M$ to $G$ is the subset of $C^\infty(P)$ whose elements solve
(\ref{t2}). Using Einstein's summation convention, we have the following coordinate expressions
\begin{equation}\label{t3}
L(j^1\phi)={\frac{1}{2}}g^{ij}{\frac {\partial \phi^{\alpha}}{\partial x^i}}{\frac {\partial \phi^{\beta}}{\partial x^j}}h_{\alpha \beta},
\end{equation}
and for (\ref{t2})
\begin{equation}\label{t4}
g^{ij} \left({\frac {\partial ^2 \phi ^{\gamma}}{\partial x^i \partial x^j}}-\Gamma
 ^k _{ij} {\frac {\partial \phi ^\gamma}{\partial x^k}}+\tilde{\Gamma} ^\gamma _{\alpha \beta}{\frac{\partial \phi ^\alpha}{\partial x^i}}{\frac{\partial \phi ^\beta}{\partial x^j}}\right)=0, \hspace{4mm} 1\leq \gamma \leq m,
\end{equation}
where $\Gamma ^k _{ij}$ , $\tilde{\Gamma} ^{\gamma} _{\alpha \beta }$ denote the Christoffel symbols of the Levi-Civita connections of $g$ and $h$.
We shall derive the reduced form of $(\ref{t2})$ for two specific
cases:  $G={\mathbb R}$, and $G={\mathbb S}^3 \cong SU(2)
\stackrel{2:1}{\cong }SO(3)$.   

For the case that $G={\mathbb R}$, the abelian group of translations, 
we choose the trivial connection for $P$.  The divergence operator div$^M$
is naturally defined by the metric $g$ and
its associated Riemannian connection.   In this case,
$(TP/G)_x\simeq T_xM \times {\mathbb R}$ and $(\text{ad}P)_x\simeq {\mathbb R}$. 
Let $\sigma = q(T\phi)$, so that $\sigma_x:T_xM \rightarrow T_xM \times
{\mathbb R}$, acting as the identity on the first factor. Then, $\sigma$ can be considered as a 1-form with local expression $\sigma =p_i\mathrm{d}x^i$, $p_i =\partial \phi /\partial x^i$.

The Lagrangian $L$ is clearly ${\mathbb R}$-invariant. 
Denoting by $l$ the projection of $L$ to ${\mathcal C}(P)$, Theorem
\ref{thm1} asserts that  $\sigma$ satisfies 
$$ \text{div}^M \frac{\delta l}{\delta\sigma}=0,$$
or, in coordinates,
$$\mathrm{div}^M (g^{kj}p_j {\frac{\partial}{\partial x^k}})={\frac{\partial (g^{ij} p_j)}{\partial x^i}}+\Gamma ^i _{ik} g^{kj}p_j=0,$$
since $l(\sigma)={\frac {1}{2}}g^{ij}p_ip_j$.
It is straightforward to check that the above equation together with vanishing curvature 
$$ {\frac{\partial p_i}{\partial x^j}}={\frac{\partial p_j}{\partial x^i}}$$
and $p_i=\partial \phi  /\partial x^i$, is equivalent to formula (\ref{t4}) for $\gamma =1$ and $\tilde{\Gamma}^{\gamma}_{\alpha \beta}=0$, as is stated in Theorem (\ref{thm2}).

For the case $G={\mathbb S}^3$, we make the identifications
$(TP/G)_x \cong T_xM \times {\mathfrak su}(2)$, 
$(\text{ad}P)_x \cong {\mathfrak su}(2)$ and ${\mathcal C}(P)\simeq T^*M\otimes 
{\mathfrak su}(2)$. Then, $\sigma =q(T\phi)$ can be considered as a 1-form taking
 values in ${\mathfrak su}(2)$. Let $\{E_1,E_2,E_3\}$ be a basis of 
${\mathfrak su}(2)$, then $\sigma$ can be written as $\sigma(x)=p^{\alpha}_i {\mathrm{d}}x^i \otimes E_{\alpha}$ with $p^{\alpha}_i\otimes E_{\alpha}=\partial \phi /\partial x^i=T\phi (\partial /\partial x^i)$.

	The Lagrangian $L$ is $SU(2)$-invariant and its projection to $\mathcal{C}(P)$ is
$$l(\sigma)={\frac{1}{2}}g^{ij}p^{\alpha}_i p^{\beta}_j h_{\alpha \beta}.$$
 Then
$${\frac{\delta l}{\delta \sigma}}=g^{ij}p^{\alpha}_i h_{\alpha \beta}{\frac{\partial}{\partial x^j}}\otimes E^{\beta},$$
and its usual divergence is
$${\mathrm{div}}{\frac{\delta l}{\delta \sigma}}=\left({\frac{\partial}{\partial x^j}}\left(g^{ij}p^{\alpha}_i h_{\alpha \beta} \right)+\Gamma ^k _{kj} g^{ij} p^{\alpha}_i h_{\alpha \beta} \right)\otimes E^{\beta}.$$
The coadjoint map can be written in coordinates as 
$$\langle {\mathrm{ad}}^*_{\sigma}{\frac{\delta l}{\delta \sigma}},E_{\beta}\rangle=\langle {\frac{\delta l}{\delta \sigma}},[\sigma,E_{\beta}]\rangle=g^{ij}p^{\alpha}_i p^{\rho}_jc^{\gamma}_{\rho \beta}h_{\alpha \gamma}.$$
Then, Euler-Poincar\'{e} equations for the trivial connection on $M\times SU(2)$ are (Theorem (\ref{thm1}))
$${\frac{\partial}{\partial x^j}}\left(g^{ij}p^{\alpha}_i h_{\alpha \beta} \right)+\Gamma ^k _{kj} g^{ij} p^{\alpha}_i h_{\alpha \beta}=-g^{ij}p^{\alpha}_i p^{\rho}_jc^{\gamma}_{\rho \beta}h_{\alpha \gamma}.$$
The above system of equations together with vanishing curvature
$${\frac{\partial p^{\gamma}_i}{\partial x^j}}-{\frac{\partial p^{\gamma}_j}{\partial x^i}}+p^{\alpha}_i p^{\beta}_j c^{\gamma}_{\alpha \beta}=0 \hspace{4mm} \forall i,j=1,...,n;\hspace{2mm}\gamma=1,2,3,$$
and $p^{\alpha}_i\otimes E_{\alpha}=\partial \phi /\partial x^i$ are equivalent to equations (\ref{t4}), as is asserted in Theorem (\ref{thm2}).

\end{document}